\documentclass[12pt, a4paper]{article}
\usepackage{graphicx, url, amsmath, amsfonts}
\pdfoutput=1

\newcommand{\PP}{{\mathbb P}}
\newcommand{\C}{{\mathbb C}}
\newcommand{\R}{{\mathbb R}}

\newcommand{\Z}{{\mathbb Z}}
\newcommand{\Q}{{\mathbb Q}}
\newcommand{\Schwarzian}{{\mathcal{S}}}
\newcommand{\HGF}{{_2F_1}}

\DeclareMathOperator{\im}{{im}}

\DeclareMathOperator{\Hes}{{\mathcal{H}}}
\DeclareMathOperator{\Jac}{{\mathcal{J}}}
\DeclareMathOperator{\Aut}{{Aut}}

\numberwithin{equation}{section}

\title{On Klein's Icosahedral Solution of the Quintic}
\author{Oliver Nash}

\begin{document}
\maketitle

\begin{abstract}
  We present an exposition of the icosahedral
  solution of the quintic equation first described in Klein's classic work
  \lq Lectures on the icosahedron and the solution of equations of the fifth degree\rq.
  Although we are heavily influenced by Klein we follow a slightly different
  approach which enables us to arrive at the solution more directly.
\end{abstract}

\section{Introduction}
In 1858, Hermite published a solution of the quintic equation
using modular functions \cite{Hermite}. His work received considerable
attention at the time and shortly afterward Kronecker \cite{Kronecker}
and Brioschi \cite{Brioschi} also
published solutions, but it was not till Klein's seminal work \cite{Klein} in 1884 that
a comprehensive study of the ideas was provided.

Although there is no modern work covering all of the material in \cite{Klein},
there are several noteworthy presentations of some of the main ideas. These include an
old classic of Dickson \cite{Dickson} as well as Slodowy's article \cite{MR879280} and the
helpful introduction he provides in the reprinted edition
\cite{KleinSlodowy} of \cite{Klein}. In addition Klein's solution is discussed in
both \cite{MR1863996}, \cite{MR1901214} as well as \cite{Hunt}. Finally the geometry is outlined
briefly in \cite{MR1471703} and a very detailed study of a slightly different approach is
presented in the book \cite{Shurman}.

Perhaps surprisingly, we believe there is room for a further exposition of the quintic's icosahedral
solution. For one thing, all existing discussions arrive at the solution of the quintic
indirectly as a result of first studying quintic resolvents of the icosahedral
field extension. Even Klein admits that he arrives at the solution~\lq somewhat
incidentally\rq\footnote{His words in the original German are \lq gewissermassen zuf\"alligerweise\rq.}
and each of the accounts listed above, except \cite{MR1471703} and \cite{Shurman}, exactly follow in
Klein's footsteps. In addition, we
believe the icosahedral solution deserves a short, self-contained account.

We thus follow Klein closely but take a direct approach to the solution of the quintic, bypassing
the study of resolvents of the icosahedral field extension. In fact our approach is
closely related to Gordon's work \cite{Gordon} and indeed Klein discusses the
connection (see \cite{Klein} part II chapter III \S 6)
but, having already achieved his goal by other means, he contents himself with an outline.

This direct approach enables us to present
the solution rather more concisely than elsewhere and we hope this may
render it more accessible; part of our motivation for writing these notes was provided by
\cite{MOModernKlein}. In addition our derivation of the icosahedral invariant
of a quintic produces a different expression than that which appears elsewhere and which is more
useful for certain purposes (for example our formula can be easily evaluated along the
Bring curve).

Finally it is worth highlighting the geometry that connects
the quintic and the icosahedron. Using a radical transformation, a
quintic can always be put in the form
$y^5 + 5\alpha y^2 + 5\beta y + \gamma = 0$. The vector of ordered roots of such a quintic lies
on the quadric surface $\sum y_i = \sum y_i^2 = 0$ in $\PP^4$ and the reduced Galois group
$A_5$ acts on the two families of lines in this doubly-ruled
surface by permuting coordinates. The $A_5$ actions on these families, parameterized by $\PP^1$,
are equivalent to the action of the group of rotations of
an icosahedron on its circumsphere and the quintic thus defines a point in the quotients ---
the icosahedral invariants of a quintic. We discuss this in detail below but first we
fix some notation and collect those facts about the icosahedron that we will need.

\noindent{\bf Acknowledgment} The author gratefully acknowledges numerous excellent suggestions
contained in an impressively thorough anonymous referee's report.

\section{The icosahedron}\label{IcosCoverSect}
Given a regular icosahedron in $\R^3$, we can identify its circumsphere $S$ with the extended complex
plane, and so also with $\PP^1$, using the usual stereographic projection:
$(x, y, z) \mapsto \frac{x + iy}{1 - z}$.
Orienting our icosahedron appropriately, the 12 vertices have complex coordinates:
\begin{align}\label{vertex_coords}
  0,~\epsilon^{\nu}(\epsilon + \epsilon^{-1}),~\infty,~\epsilon^{\nu}(\epsilon^2 + \epsilon^{-2})
  \qquad \nu=0, 1, \ldots, 4
\end{align}
where $\epsilon = e^{2\pi i/5}$.

\begin{figure}
  \begin{center}
    \includegraphics{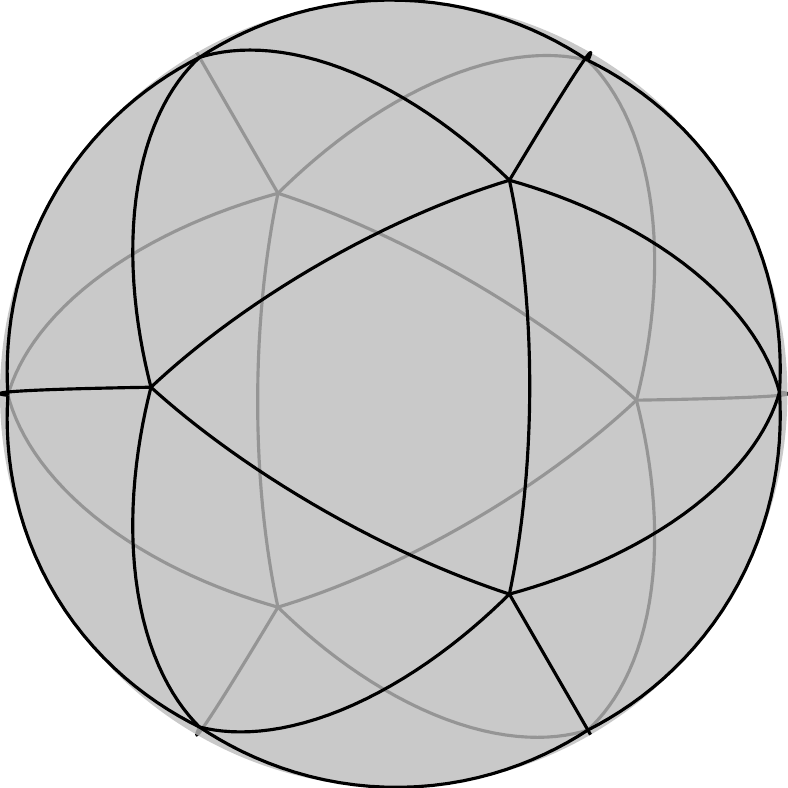}
    \caption{The icosahedron, projected radially onto its circumsphere.}
    \label{Figure1}
  \end{center}
\end{figure}

\begin{figure}
  \begin{center}
    \includegraphics{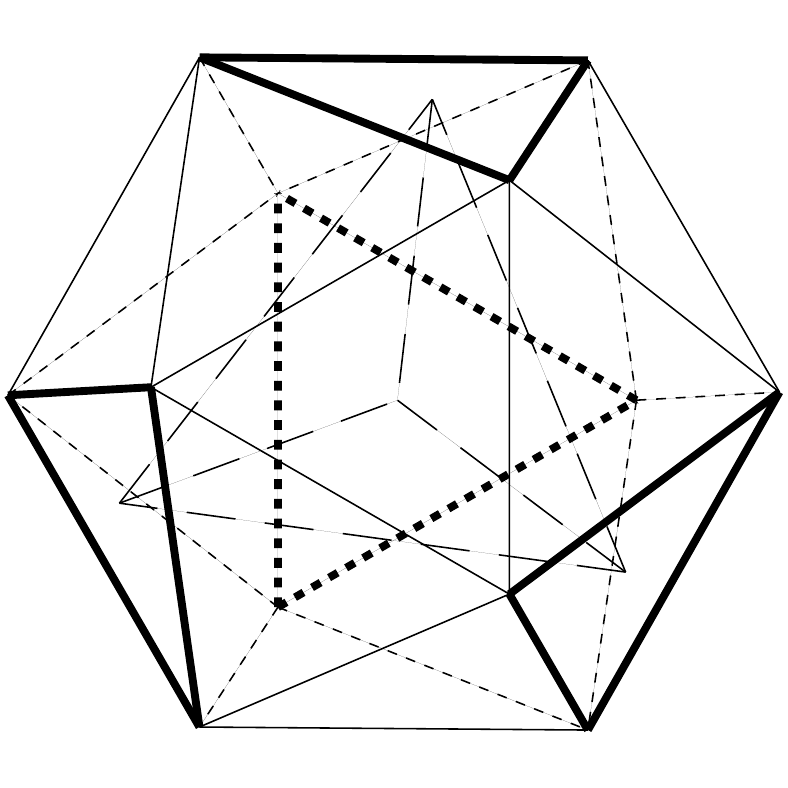}
    \caption{The icosahedron with inscribed tetrahedron.}
    \label{Figure2}
  \end{center}
\end{figure}

Projecting radially from the centre, we can
regard the edges and faces of the icosahedron as subsets of $S$.
With the sole exception of figure \ref{Figure2}, we shall always regard the faces and edges
as subsets of $S \simeq \C \cup \infty$.
The picture of the icosahedron we should have in mind is thus
similar to figure \ref{Figure1}.

We may inscribe a tetrahedron in an icosahedron by placing a tetrahedral vertex at the
centre of 4 of the 20 icosahedral faces as shown in figure \ref{Figure2}.
Note that for each icosahedral vertex, exactly one of the 5 icosahedral faces
to which it belongs has
a tetrahedral vertex at its centre. If we pick an axis joining two
antipodal icosahedral vertices, we can consider the 5 inscribed tetrahedra
obtained by rotating this configuration through $2\pi\nu/5$ for $\nu = 1, \ldots, 5$.
None of these tetrahedra have any vertices in common and so each of the 20
faces of the icosahedron are labeled by a number $\nu \in \{1, \ldots, 5\}$.
Figure 3 exhibits such a numbering after stereographic projection.
\begin{figure}
  \begin{center}
    \includegraphics{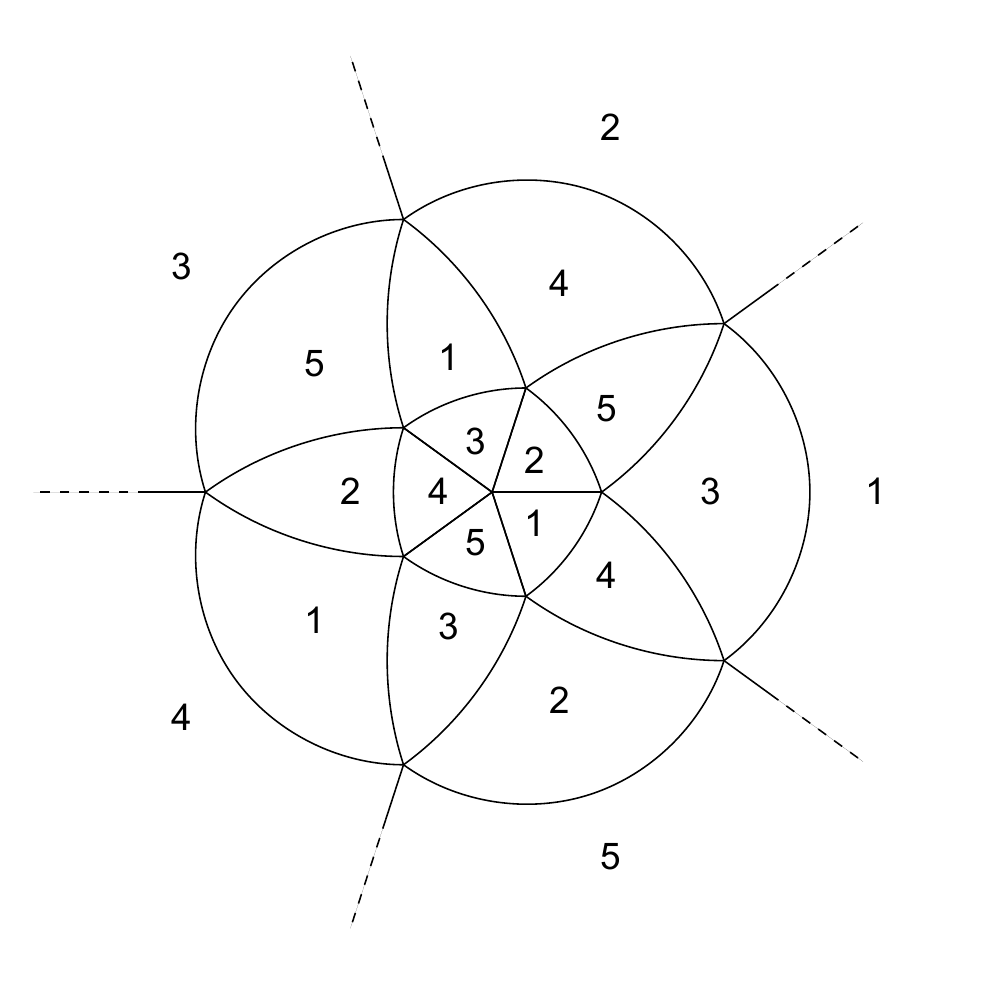}
    \caption{Tetrahedral face numbering of the icosahedron
    under stereographic projection (the outer radial lines meet at
    the vertex at infinity).}
    \label{Figure3}
  \end{center}
\end{figure}

The group $\Gamma$ of rotations of the icosahedron acts transitively on the set of 20 faces
with stabilizer of order 3 at each face and so has order 60. $\Gamma$ also
acts faithfully
on the set of 5 tetrahedra as constructed above and so we obtain
an embedding $\Gamma \hookrightarrow S_5$. Since $A_5$ is the only subgroup of $S_5$
of order 60 we must thus have:
\begin{align*}
  \Gamma \simeq A_5
\end{align*}

It will be useful later to have explicit generators for $\Gamma$. Thus let
$S$ be a rotation through $2\pi/5$ about the axis of symmetry joining the antipodal vertex
pair $0, \infty$ and let
$T$ be the rotation through $\pi$ about the axis of symmetry joining the midpoints of the
antipodal edge pair $[0, \epsilon+\epsilon^{-1}], [\infty, \epsilon^2 + \epsilon^{-2}]$.
Using the face numbering in figure \ref{Figure3}, $S, T$ correspond to the permutations:
\begin{alignat}{3}\label{ST_perms}
  S &= (12345) & \qquad T &= (12)(34)
\end{alignat}
We note in passing that since these two permutations generate $A_5$ we can use this
to see that the action on tetrahedra is faithful.

In addition, under the
embedding of symmetry groups: $\Gamma \hookrightarrow PSL(2, \C)$
associated to the identification of the circumsphere of the icosahedron with $\PP^1$
we have:
\begin{alignat}{4}\label{ST_matrices}
  S &= \left[
  \begin{array}{cc}
    \epsilon^3 & 0\\
    0          & \epsilon^2\\
  \end{array}\right]
  &\qquad&
  T &= \frac{1}{\sqrt{5}}\left[
  \begin{array}{cc}
    -(\epsilon - \epsilon^4) & \epsilon^2 - \epsilon^3\\
    \epsilon^2 - \epsilon^3  & \epsilon - \epsilon^4
  \end{array}\right]
\end{alignat}

We wish to study the branched covering:
\begin{align*}
  \PP^1 \to \PP^1 / \Gamma
\end{align*}
It will be useful to generalize slightly and work over a field that
is not necessarily $\C$. We thus take our base field $k$ to be any subfield of $\C$
containing $\epsilon = e^{2\pi i/5}$.
Note that requiring $\epsilon \in k$ avoids arithmetic issues
discussed by Serre in \cite{MR604216}\footnote{For the benefit of those consulting
\cite{MR604216}, we note that $\epsilon \in k$ guarantees $\sqrt{5} = 1 + 2(\epsilon +
\epsilon^{-1}) \in k$ and that $-1$ is a sum of two squares in $k$:
$((\epsilon - \epsilon^{-1})/\sqrt{5})^2 + ((\epsilon^2 - \epsilon^{-2})/\sqrt{5})^2 = -1$.}.

We wish to construct an explicit isomorphism $\PP^1/\Gamma \simeq \PP^1$. In general
the procedure for computing the quotient of a projective variety by a finite group
is to compute the ring of invariant elements of its homogeneous coordinate ring and then,
if necessary, replace this with a regraded subring that is generated by elements in
degree 1 (see \cite{MR1416564} for an elementary discussion).
We thus begin by computing $k[z_1, z_2]^{\Gamma}$.

Consider the possible stabilizer subgroups for the action of $\Gamma$
on $\PP^1$. The action is free
except on the three exceptional orbits that correspond to the sets of vertices, edge midpoints and
face centres where it has stabilizer
subgroups of order $5, 2, 3$ respectively. Each of these exceptional
orbits is the divisor of an invariant homogeneous polynomial.
In fact we may need to pass to an extension of $k$ in order for the
edge midpoints and face centres to be defined\footnote{The edge midpoints are
the orbit of $-i(\epsilon - \epsilon^4) + (\epsilon^2 + \epsilon^4)$
and the face centres are the orbit of $1 - \omega\epsilon - \omega^2\epsilon^4$ where
$\omega = e^{2\pi i/3}$. Thus the points of all exceptional orbits are defined iff $k$
contains a primitive $60^{\rm{th}}$ root of unity.} but this is not a problem for as
we shall see, their corresponding polynomials are defined over $k$.
Using \eqref{vertex_coords} we can calculate the polynomial corresponding
to the vertices directly. We obtain:
\begin{align}\label{invt_f_defn}
  f(z_1, z_2) &= z_1 z_2 (z_1^{10} + 11z_1^5 z_2^5 - z_2^{10})
\end{align}
To find the polynomials corresponding to the edge midpoints and face
centres, we use the Hessian and Jacobian covariants of $f$. Thus
recall (see e.g., \cite{Shurman} or \cite{MR2004511}) that
if $f, g$ are invariant polynomials in two variables, then the following
are also invariant polynomials:
\begin{alignat*}{4}
  \Hes(f)    &= \left|
                  \begin{array}{cc}
                    f_{,_{11}} & f_{,_{12}}\\
                    f_{,_{21}} & f_{,_{22}}
                  \end{array}
                \right|
  &\qquad&
  \Jac(f, g) &= \left|
                  \begin{array}{cc}
                    f_{,_1} & f_{,_2}\\
                    g_{,_1} & g_{,_2}
                  \end{array}
                \right|
\end{alignat*}
where $f_{,_i}$ denotes the partial derivative with $f$ with respect to its $i^{\rm th}$
argument and $f_{,_{ij}}$ is the iterated partial.
Clearly $\deg\Hes(f) = 2\deg f - 4$ and $\deg\Jac(f, g) = \deg f + \deg g - 2$.
Now let\footnote{We think it worthwhile following the notation of
\cite{Klein} as closely as possible to aid the reader who wishes to compare. It is
unfortunate that we must thus use $T$ to denote both the rotation mentioned above as
well as the invariant polynomial introduced here but we trust that context will protect
us from confusion.}:
\begin{alignat*}{4}
  H &= \frac{1}{121}\Hes(f) &\qquad& T &= \frac{1}{20}\Jac(f, H)
\end{alignat*}
In view of their degrees, $H, T$ must be the invariant polynomials corresponding
to the face centres and edge midpoints respectively.
Straightforward computation reveals:
\begin{align}\label{invts_TH_defn}
  \begin{split}
    H(z_1, z_2) &= -(z_1^{20} + z_2^{20}) + 228(z_1^{15} z_2^5 - z_1^5 z_2^{15}) - 494z_1^{10} z_2^{10}\\
    T(z_1, z_2) &= (z_1^{30} + z_2^{30}) + 522(z_1^{25} z_2^5 - z_1^5 z_2^{25}) -
                   10005(z_1^{20} z_2^{10} + z_1^{10} z_2^{20})
  \end{split}
\end{align}

We claim that $H^3, T^2$ form a basis for the vector space of invariant polynomials
of degree 60. It is sufficient to establish this over $\C$ since by descent:
\begin{align}\label{descent}
  \C[f, H, T] \cap k[z_1, z_2] &= k[f, H, T]
\end{align}
Firstly note that,
a non-zero, degree-60 invariant polynomial $p$ vanishes on a unique $\Gamma$ orbit.
Now consider $aH^3 + bT^2$ for scalars $a, b$, not both 0. Since the condition for this to vanish
at a given point is just a linear condition on $a, b$, we can arrange for it to vanish
at a root of $p$. By $\Gamma$-invariance it thus has the same divisor as $p$
and so must coincide with $p$ up to scale.

In particular, there must be a linear relationship $aH^3 + bT^2 = cf^5$. Evaluating
at 0 yields $a=b$ and without loss of generality we may assume $a=b=1$.
Expanding and comparing coefficients of $z_1^{60}$ we find
$c = 1728$ and thus obtain the syzygy:
\begin{align}\label{syzygy}
  H^3 + T^2 = 1728f^5
\end{align}

Now if $p \in k[z_1, z_2]^\Gamma$ is any non-zero element then, passing if
necessary to the splitting field, the divisor of $p$ is a sum of $\Gamma$ orbits,
repeated according to multiplicity. By the above, there is a linear combination
of $H^2, T^3$ vanishing on any free orbit. We can thus use $f, H, T$ to construct
an invariant polynomial with same divisor as $p$ and so obtain:
\begin{align*}
  p &= c f^{e_1}H^{e_2}T^{e_3} \prod_j (a_j H^2 + b_j T^3)
\end{align*}
for scalars $c, a_j, b_j$ and natural numbers $e_i$. As before, by \eqref{descent}
we thus have:
\begin{align*}
  k[z_1, z_2]^\Gamma = k[f, H, T]
\end{align*}

Thus if we define the graded ring $A$ as:
\begin{align*}
  A = k[x, y, z] / (1728 x^5 - y^2 - z^3)
\end{align*}
where $x, y, z$ are pure transcendental and are given gradings of $12, 30, 20$
respectively, then we have a natural surjection of graded rings:
\begin{align*}
   A \to k[z_1, z_2]^\Gamma
\end{align*}
Since any surjection between finite-dimensional integral domains of the same
dimension is necessarily an isomorphism (pull back a maximum-length chain of prime ideals)
this map must in fact be an isomorphism.

Finally, note that if $A = \oplus_{n\ge 0}A_n$ is the graded decomposition of $A$
then $A$ is not generated in degree 1 (indeed $A_1 = 0$) so we pass to:
\begin{align*}
  A^{(60)} = \bigoplus_{n\ge 0} A_{60n}
\end{align*}
where we define the grading as $A^{(60)}_n = A_{60n}$.
Then $A^{(60)}$ is generated in degree 1. Thus we can take $A^{(60)}$ as the homogeneous
coordinate ring of $\PP^1/\Gamma$.

In fact $A^{(60)} \simeq k[x^5, y^2]$
is a polynomial algebra and so $\PP^1/\Gamma \simeq \PP^1$.
This shows that the map:
\begin{align*}
  \PP^1 &\to \PP^2\\
  [z_1, z_2] &\mapsto [H^3, T^2, 1728f^5]
\end{align*}
to the line $\{a + b = c\} \subset \PP^2$ is a quotient map for the action of $\Gamma$.
Following Klein, we identify this line
with $\PP^1$ by sending $[0, 1, 1], [1, 0, 1], [1, -1, 0]$ to $0, 1, \infty$
respectively. This realizes the quotient map as:
\begin{align}\label{icos_quot_map_defn}
  I &= \frac{H^3}{1728 f^5}
\end{align}

If we were now to follow the usual approach to the icosahedral solution of the quintic,
we would next study quintic resolvents,
of the Galois extension $k(\PP^1) \supset k(\PP^1)^\Gamma$. These are obtained by taking
index-5 subgroups of the Galois group $A_5$ corresponding to the tetrahedron.
However, as we have said, we follow a slightly different
approach and so immediately turn our attention to the solution of the quintic.

\section{Tschirnhaus and the canonical equation}\label{Tschirnhaus_sect}
A common approach when solving the polynomial equation:
\begin{align}\label{gen_poly_eqn}
  x^n + a_1 x^{n-1} + \cdots + a_n &= 0
\end{align}
is to begin by making the affine substitution $y = x + a_1/n$ and so eliminate
the term of degree $n-1$. This substitution is a special case of the so-called Tschirnhaus
transformation \cite{Tschirnhaus} in which $y$ is allowed to be
a polynomial expression $q$ in $x$. If $\alpha_i$ are the roots of \eqref{gen_poly_eqn},
the coefficients of the transformed equation:
\begin{align*}
  \prod_i \left(y-q(\alpha_i)\right) &= 0
\end{align*}
are polynomials in the $a_i$ by $S_n$-invariance (or Newton's identities).

Using a Tschirnhaus transformation we can simultaneously eliminate
further terms in the original polynomial. For example if $n \ge 3$ and $a_1 = 0$,
it is easy to check that the substitution:
\begin{align*}
  y = x^2 + b_1 x + b_2
\end{align*}
simultaneously eliminates the terms of degree $n-1$ and $n-2$ provided the
coefficients $b_1, b_2$ satisfy the auxiliary polynomial conditions \cite{Shurman}:
\begin{align*}
  b_2 - p_2 / n &= 0\\
  p_2 b_1^2 + 2 p_3 b_1 + (p_4 - p_2^2 / n) &= 0
\end{align*}
where $p_j = \sum_i \alpha_i^j$ are the power sums of the roots.

Thus, provided we are willing to allow ourselves the auxiliary square root necessary
to solve the above quadratic for $b_1$, we may take the general form of the
quintic to be\footnote{We include the factors of 5 for consistency with \cite{Klein}.}:
\begin{align}\label{CanonicalQuintic}
  y^5 + 5\alpha y^2 + 5\beta y + \gamma = 0
\end{align}

In fact it is possible to simultaneously
eliminate the terms of degrees $n-1$, $n-2$ and $n-3$ (where the coefficients
of the substitution are determined by polynomials of degree strictly less than $n$). Thus, as
first shown by Bring \cite{Bring} and subsequently by Jerrard \cite{Jerrard}, the
general quintic can be reduced to the so-called Bring-Jerrard form:
\begin{align*}
  y^5 + y + \gamma &= 0
\end{align*}
However it is not in this form that the quintic most easily reveals its icosahedral
connections and so, except for section \ref{BringCurveSection} and appendix
\ref{eisenstein_appendix}, we shall take the quintic in the form \eqref{CanonicalQuintic}.

\section{The icosahedral invariant}\label{IcosInvtSect}
The key to Klein's solution of the quintic is his icosahedral invariant. Working
over $\C$, we sketch the geometric idea before turning to the algebraic job
of calculating the invariant when we shall be more precise.

Thus consider the quintic \eqref{CanonicalQuintic} for $\alpha, \beta, \gamma \in \C$
not all 0. Given an ordering, we may regard the roots as homogeneous coordinates of a point
in $\PP^4$. If, in addition to $\alpha, \beta, \gamma$, we also supply
a distinguished square root of the discriminant (which we assume is non-zero):
\begin{align}\label{discrim_defn}
  D &= 3125\prod_{i<j}(y_i-y_j)^2\notag\\
    &= 108\alpha^5\gamma - 135\alpha^4\beta^2 + 90\alpha^2\beta\gamma^2 -
       320\alpha\beta^3\gamma + 256\beta^5 + \gamma^4
\end{align}
then the roots are ordered up to even permutation and so we obtain an $A_5$ orbit
in $\PP^4$ where $S_5$ acts by permuting coordinates. Furthermore,
because the quintic lacks terms of degree 4 and 3, this orbit lies in the
non-singular $S_5$-invariant quadric surface:
\begin{align*}
  Q &= \left\{\left[y_0, \ldots, y_4\right]\in\PP^4\quad|\quad\sum y_i = \sum y_i^2 = 0\right\}
\end{align*}
Now $Q$ is a doubly-ruled surface $Q \simeq \PP^1 \times \PP^1$. The $A_5$ action
sends lines to lines and so the $\PP^1$s
appearing in the double ruling come with $A_5$ actions and the ruling is equivariant.
Projection onto either factor of $Q \simeq \PP^1\times \PP^1$ and taking
$A_5$ quotient yields an invariant
in $\PP^1/A_5 \simeq \PP^1$.

We turn our attention to the calculation of these invariants. We fix
our base field as $k = \Q(\epsilon)$ and, for simplicity, we assume
that $\alpha, \beta, \gamma \in \C$ are algebraically independent over $k$.
We also let $\nabla \in \C$ be a square root of the discriminant \eqref{discrim_defn}
and let the roots of \eqref{CanonicalQuintic} be $y_1, \ldots, y_5 \in \C$. We
have the following diagram of $k$-algebra isomorphisms:
\begin{center}
\begin{tabular}{ccccc}
  $k_{hom}[Q]$ &
  $=$ &
  $\frac{k[\hat y_1, \ldots, \hat y_5]}{(\sum \hat y_i, \sum \hat y_i^2)}$ &
  $\simeq$ &
  $k[y_1, \ldots, y_5]$\\
  \\
  $\bigcup$ && $\bigcup$ && $\bigcup$\\
  \\
  $k_{hom}[Q]^{A_5}$ &
  $\simeq$ &
  $\frac{k[\hat \alpha, \hat \beta, \hat \gamma, \hat \nabla]}{(\hat\nabla^2 - D)}$ &
  $\simeq$ &
  $k[\alpha, \beta, \gamma, \nabla]$\\
  \\
  $\bigcup$ && $\bigcup$ && $\bigcup$\\
  \\
  $k_{hom}[Q]^{S_5}$ &
  $\simeq$ &
  $k[\hat \alpha, \hat \beta, \hat \gamma]$ &
  $\simeq$ &
  $k[\alpha, \beta, \gamma]$
\end{tabular}
\end{center}
Here $k_{hom}[Q]$ is the homogeneous coordinate ring of $Q$,
the $\hat y_i$ are pure transcendental, $\hat\alpha, \hat\beta, \hat\gamma$
are the elementary symmetric functions in the $\hat y_i$ of degrees 3, 4, 5 respectively,
$\hat\nabla$ is pure transcendental and $D$ is the discriminant polynomial \eqref{discrim_defn}
in the $\hat{}$ variables. The maps are those suggested by the notation (i.e., remove
$\hat{}$ s). It is straightforward to verify the various maps are isomorphisms
using the following well-known facts:
\begin{itemize}
  \item The elementary symmetric functions are algebraically independent.
  \item If $e_1, e_2$ are the elementary symmetric functions in the variables
        $\hat y_1, \ldots \hat y_5$ of degrees 1, 2 then:
        $k[\hat y_1, \ldots, \hat y_5]^{S_5} \simeq
        k[e_1, e_2, \hat\alpha, \hat\beta, \hat\gamma]$
        and $k[\hat y_1, \ldots, \hat y_5]^{A_5} \simeq
        k[e_1, e_2, \hat\alpha, \hat\beta, \hat\gamma, \hat\nabla]/(\hat\nabla^2 - D)$.
  \item If a $k$-algebra $A$ carrying an action of a finite group $G$ has $G$-invariant
        ideal $\mathfrak{a}$, then $(A/\mathfrak{a})^G \simeq A^G/\mathfrak{a}^G$.
	(Indeed, $1/|G|\sum_{g\in G} g\cdot a$ is a lift of any $[a] \in (A/\mathfrak{a})^G$.)
  \item A surjection between integral domains of the same dimension is an isomorphism.
\end{itemize}

To define the icosahedral invariant in this setting, we need the algebraic
expression of the double ruling. To this end we introduce:
\begin{align}\label{p_basis_defn}
  p_k = \sum_j \epsilon^{kj}y_j
\end{align}
and, letting $k[u, v]_n$ denote the degree-$n$ component of the graded
ring $k[u, v]$, we define an isomorphism of graded $k$-algebras:
\begin{align}\label{Segre_product}
  k[y_1, \ldots, y_5] &\to \bigoplus_{n\ge 0} k[\lambda_1, \lambda_2]_n \otimes k[\mu_1, \mu_2]_n
\end{align}
by making the identifications:
\begin{align}\label{lambda_mu_defn}
  p_1 =  5\lambda_1\mu_1\quad p_2 = -5\lambda_2\mu_1\quad p_3 =  5\lambda_1\mu_2\quad p_4 =  5\lambda_2\mu_2
\end{align}
Furthermore, a computation reveals that if we let $A_5$ act on $k[\lambda_1, \lambda_2]$ using
the formulae \eqref{ST_matrices} and act on $k[\mu_1, \mu_2]$ using the same
formulae after replacing $\epsilon$ with $\epsilon^2$, then \eqref{Segre_product}
is an $A_5$-equivariant isomorphism.
Finally we use \eqref{Segre_product} to define a full $S_5$ action on
on $\oplus_{n\ge 0}k[\lambda_1, \lambda_2]_n \otimes k[\mu_1, \mu_2]_n$. It is
sufficient to define the action of any odd permutation and a quick
computation reveals that action of $R = (1243)$ can be described by:
\begin{align}\label{odd_perm_R_action}
  ([\lambda_1, \lambda_2], [\mu_1, \mu_2]) \mapsto ([\mu_2, -\mu_1], [\lambda_1, \lambda_2])
\end{align}

We already have the algebraic expression of $\PP^1 / A_5 \simeq \PP^1$;
it is the rational map \eqref{icos_quot_map_defn}.
Thus, recalling our formulae for $f, H, T$ in section \ref{IcosCoverSect}, we define:
\begin{align*}
  f_1 = f(\lambda_1, \lambda_2)\quad f_2 = f(\mu_1, \mu_2)
\end{align*}
and similarly we define $H_1, H_2$ and $T_1, T_2$. Finally we can define the icosahedral
invariants:
\begin{align*}
  Z_i &= \frac{H_i^3}{1728f_i^5}
\end{align*}
A priori we have $Z_1 \in k(\lambda_1, \lambda_2)$ and $Z_2 \in k(\mu_1, \mu_2)$, however writing:
\begin{align}\label{Z1_formula}
  Z_1 &= \frac{H_1^3f_2^5}{1728(f_1f_2)^5}
\end{align}
and using the isomorphism \eqref{Segre_product} we may regard $Z_1$ (and similarly $Z_2$)
as an element of the splitting field $k(y_1, \ldots, y_5)$. Then by $A_5$ invariance
we have:
\begin{align*}
  Z_i \in k(y_1, \ldots, y_5)^{A_5} = k(\alpha, \beta, \gamma, \nabla)
\end{align*}

Our goal now is to compute $Z_i$ in terms of $\alpha, \beta, \gamma, \nabla$. We
deal with the numerator and denominator of \eqref{Z1_formula} separately. They each
lie in $k[\alpha, \beta, \gamma, \nabla]$ and any element
$h \in k[\alpha, \beta, \gamma, \nabla]$ can be written as:
\begin{align*}
  h &= h_s + h_a \nabla
\end{align*}
for unique polynomials $h_s, h_a \in k[\alpha, \beta, \gamma] = k[y_1, \ldots, y_5]^{S_5}$
determined by:
\begin{align}\label{symm_anti_symm_decomp}
  \begin{split}
    h_s &= (h + h^*)/2\\
    h_a\nabla &= (h - h^*)/2
  \end{split}
\end{align}
where $h^*$ is the polynomial obtained by acting on $h$ with any odd permutation.

Note that by \eqref{odd_perm_R_action} $R$ interchanges $f_1, f_2$ and so
$f_1f_2 \in k[\alpha, \beta, \gamma]$.
Since $f_1f_2$ is of degree 12, it must be a linear combination of
$\alpha^4, \beta^3, \alpha\beta\gamma$. To fix the coefficients we compare leading
coefficients as polynomials in $\lambda_i, \mu_i$. It is straightforward to verify
that:
\begin{align*}
  \alpha =& -\lambda_1^3\mu_1^2\mu_2 - \lambda_1^2\lambda_2\mu_2^3
            -\lambda_1\lambda_2^2\mu_1^3 + \lambda_2^3\mu_1\mu_2^2\\
  \beta  =& -\lambda_1^4\mu_1\mu_2^3 + \lambda_1^3\lambda_2\mu_1^4
            + 3\lambda_1^2\lambda_2^2\mu_1^2\mu_2^2 - \lambda_1\lambda_2^3\mu_2^4 +
            \lambda_2^4\mu_1^3\mu_2\\
  \gamma =& -\lambda_1^5(\mu_1^5 + \mu_2^5) + 10\lambda_1^4\lambda_2\mu_1^3\mu_2^2 -
            10\lambda_1^3\lambda_2^2\mu_1\mu_2^4 -\\
          &~10\lambda_1^2\lambda_2^3\mu_1^4\mu_2 -
            10\lambda_1\lambda_2^4\mu_1^2\mu_2^3 + \lambda_2^5(\mu_1^5 - \mu_2^5)
\end{align*}
The coefficient of $\lambda_1^{12}$ in $f_1f_2$ is 0 whereas the same coefficients
in $\alpha^4, \beta^3, \alpha\beta\gamma$
are $\mu_1^8\mu_2^4, -\mu_1^3\mu_2^9, -\mu_1^8\mu_2^4-\mu_1^3\mu_2^9$ respectively. From this
we see that we must have $f_1f_2 = A(\alpha^4 - \beta^3 + \alpha\beta\gamma)$ for some constant
$A$. Furthermore, upon noting that the coefficient of $\lambda_1^{11}\lambda_2\mu_1^{11}\mu_2$ in $f_1f_2$
is 1 whereas it is 0 in $\alpha^4, \beta^3$ and 1 in $\alpha\beta\gamma$ we learn that $A = 1$.
In other words we obtain:
\begin{align*}
  f_1 f_2  =&~\alpha^4 - \beta^3 + \alpha\beta\gamma
\end{align*}
This deals with the denominator in \eqref{Z1_formula}; we turn our attention to the numerator.

Decomposing the numerator of \eqref{Z1_formula} using \eqref{symm_anti_symm_decomp} and
recalling that our odd permutation $R$ interchanges the $f_1, f_2$ as well as $H_1, H_2$, we get:
\begin{align}\label{H13f25_pq_form}
  H_1^3f_2^5 &= \frac{H_1^3f_2^5 + H_2^3f_1^5}{2} + \frac{H_1^3f_2^5 - H_2^3f_1^5}{2}\notag\\
             &= p + \nabla q
\end{align}
where $p, q$ are polynomials in $\alpha, \beta, \gamma$.

We could now attempt to calculate $p, q$ in the same way that we calculated $f_1f_2$
above but this would be a long calculation since $p, q$
have degrees 60, 50 respectively. Instead, recall that we have the syzygies:
\begin{align*}
  T_i^2 &= 12^3f_i^5 - H_i^3
\end{align*}
Multiplying these together and rearranging we obtain:
\begin{align*}
  2p &= H_1^3 f_2^5 + H_2^3 f_1^5 = 12^3(f_1 f_2)^5 + 12^{-3}(H_1 H_2)^3 - 12^{-3} (T_1 T_2)^2
\end{align*}
We will thus have the required expression for $p$ in terms of $\alpha, \beta, \gamma$
as soon as we express $H_1H_2$ and $T_1T_2$ in these terms. To do this we use the
same procedure that we used to find $f_1f_2$ above and (admittedly with somewhat more
effort) we obtain:
\begin{align*}
  H_1 H_2  =&~\gamma^4 + 40\alpha^2\beta\gamma^2 - 192\alpha^5\gamma - 120\alpha\beta^3\gamma +
             640\alpha^4\beta^2 - 144\beta^5\\
  T_1 T_2  =&~\gamma^6 + 60\alpha^2\beta\gamma^4 + 576\alpha^5\gamma^3 - 180\alpha\beta^3\gamma^3 +
              648\beta^5\gamma^2 - 2760\alpha^4\beta^2\gamma^2 +\\
            &~7200\alpha^7\beta\gamma - 1728\alpha^{10} +
              9360\alpha^3\beta^4\gamma - 2080\alpha^6\beta^3 - 16200\alpha^2\beta^6
\end{align*}

It remains only to calculate $q$. This time the trick we use is to note that as
well as \eqref{H13f25_pq_form} above, we have $H_2^3f_1^5 = p - \nabla q$ and so:
\begin{align*}
  (H_1H_2)^3 (f_1f_2)^5 &= p^2 - \nabla^2 q^2
\end{align*}
It follows that taking our above polynomial expressions for
$H_1H_2, f_1f_2, \nabla^2 = D(\alpha, \beta, \gamma), p$
we must find a factorization of $((H_1H_2)^3(f_1f_2)^5 - p^2) / D(\alpha, \beta, \gamma)$.
From this we determine:
\begin{align*}
  2q =\pm&~(-8\alpha^5\gamma - 40\alpha^4\beta^2 + 10\alpha^2\beta\gamma^2 + 45\alpha\beta^3\gamma
        - 81\beta^5 - \gamma^4) \cdot\\
      &~(64\alpha^{10} + 40\alpha^7\beta\gamma - 160\alpha^6\beta^3 + \alpha^5\gamma^3 -\\
      &~5\alpha^4\beta^2\gamma^2 + 5\alpha^3\beta^4\gamma - 25\alpha^2\beta^6 - \beta^5\gamma^2)
\end{align*}
The two signs corresponding to the two invariants: $Z_1, Z_2$.
With this formula in hand we have achieved our goal of expressing $Z_i$ in terms of
$\alpha, \beta, \gamma, \nabla$.

\section{Obtaining the roots}
Given a quintic \eqref{CanonicalQuintic} with icosahedral invariant $Z = Z_1$, we
know that the roots of the degree-60 polynomial equation in $z$ over
$k(\alpha, \beta, \gamma, \nabla)$:
\begin{align}\label{icos_poly_eqn}
  H(z, 1)^3 - 1728Z f(z, 1)^5 = 0
\end{align}
all lie in the splitting field of the quintic. Indeed $z = \lambda_1/\lambda_2 = p_3/p_4$
is a root and all others are
obtained by the action of the Galois group. In the next section we will show how
to obtain a root of \eqref{icos_poly_eqn}. Here we show how $z$ enables
us to find the roots of our quintic equation using only rational expressions.

Thus note that by \eqref{p_basis_defn}, \eqref{lambda_mu_defn} we have:
\begin{align}\label{y_from_lambda_mu}
  y_\nu &= \epsilon^{4\nu}\lambda_1\mu_1 - \epsilon^{3\nu}\lambda_2\mu_1 +
           \epsilon^{2\nu}\lambda_1\mu_2 + \epsilon^\nu \lambda_2\mu_2
\end{align}
We now take up an idea of Gordon's \cite{Gordon} and note that if we can find $A_5$-invariant
forms that are linear in
$\mu_i$ then we can use these to eliminate the $\mu_i$ in \eqref{y_from_lambda_mu}
and so express $y_\nu$ in terms of just
$\alpha, \beta, \gamma, \lambda_1, \lambda_2$. To do this we
enlarge the ring of invariant polynomials we are studying from
$\oplus_{n\ge 0}k[\lambda_1, \lambda_2]_n \otimes k[\mu_1, \mu_2]_n$
to the full tensor product $k[\lambda_1, \lambda_2, \mu_1, \mu_2]$.
The two invariant forms linear in $\mu_i$ of lowest degree in $\lambda_i$ are:
\begin{align}\label{MN_defn}
  \begin{split}
  N_1 &= (7\lambda_1^5 \lambda_2^2 + \lambda_2^7)\mu_1 + (-\lambda_1^7 + 7\lambda_1^2\lambda_2^5)\mu_2\\
  M_1 &= (\lambda_1^{13} - 39\lambda_1^8\lambda_2^5 - 26\lambda_1^3\lambda_2^{10})\mu_1 +
         (-26\lambda_1^{10}\lambda_2^3 + 39\lambda_1^5\lambda_2^8 + \lambda_2^{13})\mu_2
  \end{split}
\end{align}
There are a number of ways to derive these expressions. We
follow Gordon \cite{Gordon} and use transvectants. We thus recall (see for example
\cite{CrawleyBoevey} or \cite{MR2004511}) that if $f, g$
are two homogeneous polynomials in $\lambda_1, \lambda_2$ then the $r^{\rm th}$ transvectant
of $f, g$ is given by:
\begin{align*}
  (f, g)_r = \sum_{i=0}^r \frac{(-1)^i}{i!(r-i)!}
                          \frac{\partial^r f}{\partial \lambda_1^{r-i}\partial \lambda_2^i}
                          \frac{\partial^r g}{\partial \lambda_1^i\partial \lambda_2^{r-i}}
\end{align*}
We extend this to homogeneous polynomials in both $\lambda_i, \mu_i$ using bilinearity, i.e.,
if:
\begin{align*}
  f &= \sum_{i,j} f_{ij} \lambda_1^i\lambda_2^{a-i}\mu_1^j\mu_2^{b-j}\\
  g &= \sum_{k,l} g_{kl} \lambda_1^k\lambda_2^{c-k}\mu_1^l\mu_2^{d-l}
\end{align*}
then we define the $(r, s)$-transvectant:
\begin{align*}
  (f, g)_{r,s} = \sum_{i,j,k,l} f_{ij}g_{kl}
                 (\lambda_1^i\lambda_2^{a-i}, \lambda_1^k\lambda_2^{c-k})_r
                 (\mu_1^j\mu_2^{b-j}, \mu_1^l\mu_2^{d-l})_s
\end{align*}
It is then straightforward to verify that:
\begin{align*}
  (\alpha, \beta)_{0, 3} &= 6N_1\\
  ((\alpha, \alpha)_{0, 2}, N_1)_{0, 1} &= 8M_1
\end{align*}

Note that geometrically, $N_1, M_1$ are $A_5$-equivariant branched
covers: $\PP^1 \to \PP^1$. It should be possible to exploit this
point of view to obtain an alternate derivation of $N_1, M_1$.
(E.g., since the branch locus must be $A_5$ invariant the
Riemann-Hurwitz relation greatly restricts the possible degrees.)

Returning to the task at hand we solve the $2\times 2$ system \eqref{MN_defn}
and express $\mu_i$ in terms of $M_1, N_1$ and using \eqref{y_from_lambda_mu} obtain:
\begin{align}\label{roots_in_mn}
  y_\nu = H_1^{-1} b_\nu M_1 + H_1^{-1} c_\nu N_1
\end{align}
Here $H_1$ appears as it is the determinant of the matrix which we invert and the coefficients
$b_\nu, c_\nu$ are given by:
\begin{align*}
  \left[
  \begin{array}{cc}
    b_\nu & c_\nu
  \end{array}
  \right] &=\\
  \left[
  \begin{array}{cc}
    \epsilon^{4\nu}\lambda_1 - \epsilon^{3\nu}\lambda_2 &
    \epsilon^{2\nu}\lambda_1 + \epsilon^\nu\lambda_2\\
  \end{array}
  \right]&
  \left[
  \begin{array}{cc}
    -\lambda_1^7 + 7\lambda_1^2\lambda_2^5 &
    26\lambda_1^{10}\lambda_2^3 - 39\lambda_1^5\lambda_2^8 - \lambda_2^{13}\\
    -7\lambda_1^5\lambda_2^2 - \lambda_2^7 &
    \lambda_1^{13} - 39\lambda_1^8\lambda_2^5 - 26\lambda_1^3\lambda_2^{10}\\
  \end{array}
  \right]
\end{align*}

We wish to express everything in \eqref{roots_in_mn} in terms of
$\alpha, \beta, \gamma, \nabla, \lambda_1, \lambda_2$.
We thus rewrite it so that all forms appearing have the same degree in $\lambda_i, \mu_i$:
\begin{align}\label{roots_in_bc_first}
  y_\nu = \frac{b_\nu f_1}{H_1} \cdot \frac{M_1 f_2}{f_1 f_2} +
          \frac{c_\nu T_1}{H_1 f_1^2} \cdot \frac{N_1 f_1^2 T_2}{T_1 T_2}
\end{align}

The methods described in section \ref{IcosInvtSect} then allow us to calculate:
\begin{align*}
  M_1 f_2       &= (11\alpha^3\beta + 2\beta^2\gamma - \alpha\gamma^2)/2 - \nabla\alpha/2\\
  N_1 f_1^2 T_2 &= r + \nabla s
\end{align*}
where:
\begin{align*}
  2r =&~\alpha^2\gamma^5 - \alpha\beta^2\gamma^4 + 53\alpha^4\beta\gamma^3 + 64\alpha^7\gamma^2 -
        7\beta^4\gamma^3 - 225\alpha^3\beta^3\gamma^2 -\\
      &~12\alpha^6\beta^2\gamma + 216\alpha^9\beta + 717\alpha^2\beta^5\gamma - 464\alpha^5\beta^4 -
       720\alpha\beta^7\\
  2s =&~-\alpha^2\gamma^3 + 3\alpha\beta^2\gamma^2 - 9\beta^4\gamma - 4\alpha^4\beta\gamma -
        8\alpha^7 - 80\alpha^3\beta^3
\end{align*}
and since we already have formulae for $f_1f_2$ and $T_1T_2$ we have the required
expression for $y_\nu$ in terms of $\alpha, \beta, \gamma, \nabla, \lambda_1, \lambda_2$.

Finally, it is possible to further simplify since:
\begin{align*}
  b_\nu &= \epsilon^\nu B(\epsilon^\nu\lambda_1, \lambda_2)\\
  c_\nu &= \epsilon^{3\nu} C(\epsilon^\nu \lambda_1, \lambda_2)
\end{align*}
where $B, C$ are the polynomials defined by:
\begin{align*}
  B(z_1, z_2) &= -z_2^8 - z_1z_2^7 - 7(z_1^2z_2^6 - z_1^3z_2^5 + z_1^5z_2^3 + z_1^6z_2^2)
                 + z_1^7z_2 - z_1^8\\
  C(z_1, z_2) &= B(z_1, z_2) D(z_1, z_2)\\
  D(z_1, z_2) &= -z_1^6 - 2z_1^5z_2 + 5z_1^4z_2^2 + 5z_1^2z_2^4 + 2z_1z_2^5 - z_2^6
\end{align*}
Bearing in mind that $H(\epsilon \lambda_1, \lambda_2) = H(\lambda_1, \lambda_2),
T(\epsilon \lambda_1, \lambda_2) = T(\lambda_1, \lambda_2)$ whereas
$f(\epsilon\lambda_1, \lambda_2) = \epsilon f(\lambda_1, \lambda_2)$,
we may thus rewrite \eqref{roots_in_bc_first} as:
\begin{align*}
   y_\nu = \left.\frac{B_1 f_1}{H_1}\right|_\nu \cdot \frac{M_1 f_2}{f_1 f_2} +
           \left.\frac{B_1 D_1 T_1}{H_1 f_1^2}\right|_\nu \cdot \frac{N_1 f_1^2 T_2}{T_1 T_2}
           \qquad\nu = 0, 1, \ldots, 4
\end{align*}
where $B_1 = B(\lambda_1, \lambda_2), D_1 = D(\lambda_1, \lambda_2)$ and
the notation involving $\nu$ on the right means we evaluate at $(\epsilon^\nu \lambda_1, \lambda_2)$.

In fact, $H$ contains $B$ as a factor. Thus the two occurrences of $B_1/H_1$
in the above expression can be simplified to polynomials of degree 12. There is
a geometric explanation for this: the roots of $B, D$ are, respectively, the locations of vertices
and face centres of a regular cube and the vertices of this cube are the vertices
of an inscribed tetrahedron (as shown in figure 2) together with the vertices of
its dual tetrahedron\footnote{The face centres of the cube lie at the
midpoints of 6 of the 30 icosahedron edges and so $D$ is a factor of $T$, though we make
no use of this.}.

\section{Solving the icosahedral equation}\label{icos_inv_sect}
In this section, we work over $\C$ since we need to use analytic methods.
We wish to invert the equation:
\begin{align}\label{I_inverse_problem}
  I(z) = Z
\end{align}
where $I$ is the icosahedral function \eqref{icos_quot_map_defn}.
(In this section we regard $I$ as a function of the single variable $z=z_1/z_2$.)
This problem was essentially solved by Schwarz in his 1873 paper \cite{Schwarz} where he
determined the list
parameters for which the hypergeometric differential equation has finite monodromy.
Recall that the Schwarzian derivative of an analytic function $s$ of one variable is:
\begin{align*}
  \Schwarzian s &= \left(\frac{s''}{s'}\right)' - \frac{1}{2}\left(\frac{s''}{s'}\right)^2
\end{align*}
Now $\Schwarzian s$ is invariant under M\"obius transformation (indeed this
can be used to define $\Schwarzian$) and since any two branches of a local inverse
to \eqref{I_inverse_problem} are related by a M\"obius transformation, the
Schwarzian derivative is independent of the branch. Following \cite{MR0045823}
we show how to compute $\Schwarzian s$ for a local inverse $s$ of
\eqref{I_inverse_problem}. This yields a differential equation for $s$ which we
then solve in terms of hypergeometric series.

\begin{figure}
  \begin{center}
    \includegraphics{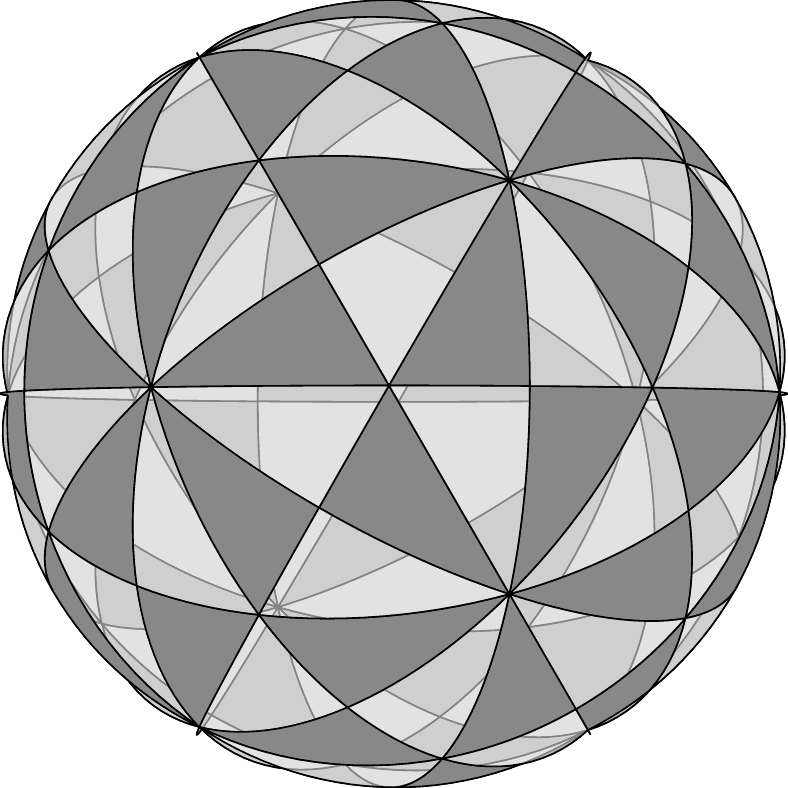}
    \caption{Icosahedral tiling of sphere.
    $I$ maps the interior of each light and dark triangle biholomorphically onto the upper
    and lower half-planes respectively.}
    \label{Figure4}
  \end{center}
\end{figure}

We begin by identifying domains of injectivity for $I$, i.e.,
fundamental domains for the action of $\Gamma$ on $\PP^1$. We thus note that
if $r$ is any reflection about a plane of symmetry of the icosahedron then since
$r$ is conjugate to $z \mapsto \bar z$ by a rotation, we must have:
\begin{align*}
  I\circ r = \bar I
\end{align*}
Since there is a plane
of symmetry through any edge of the icosahedron as well as a plane of symmetry through each
of the altitudes of any face of the icosahedron, it follows that the edges and altitudes
of the faces of the icosahedron constitute the preimage of $\R\PP^1 = \R\cup\infty$ under $I$.
The altitudes divide each face into six spherical triangles with angles $\pi / \nu_i$ where:
\begin{align*}
  \nu_1 = 2\quad\nu_2 = 3\quad\nu_3 = 5
\end{align*}
$I$ sends the vertices of each triangle to $0, 1, \infty$ (indeed we used this property
to specify $I$) and is injective on the interior. It maps three of them
biholomorphically to upper half space $H^+$ and three of them biholomorphically to
lower half space $H^-$, according to whether their vertices are sent to $0, 1, \infty$
in anti-clockwise or clockwise order respectively. Subdividing faces like this,
figure \ref{Figure1} becomes figure \ref{Figure4}.

\begin{figure}
  \begin{center}
    \includegraphics{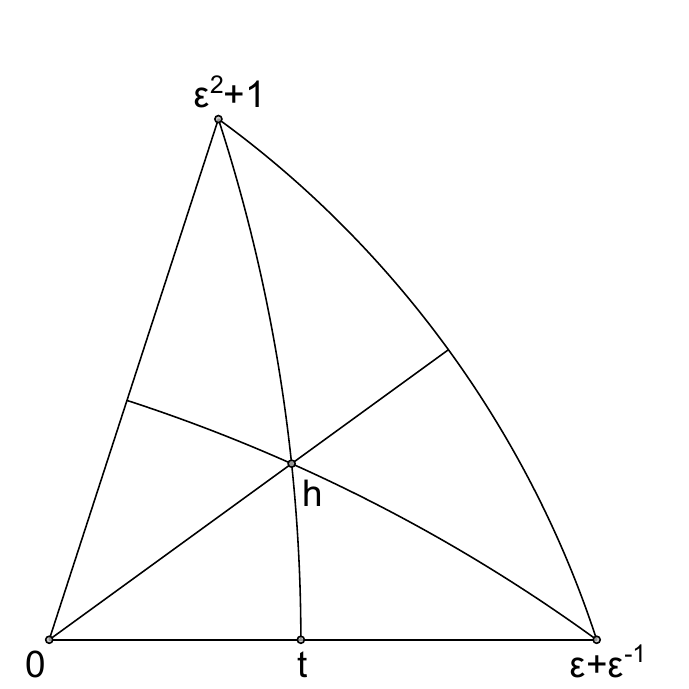}
    \caption{Domains of injectivity for $I$ under stereographic projection. The
    points $h, t$ are the images of the face centre and edge midpoint respectively.}
    \label{Figure5}
  \end{center}
\end{figure}

The subdivision of the face with vertices
$0, \epsilon + \epsilon^{-1}, \epsilon^2 + 1$ under stereographic projection
is shown in figure \ref{Figure5}.
We shall construct an inverse for the restriction of $I$ to the interior of the
triangle $\mathcal{T}$ with vertices $0, t, h$ (in the notation of figure \ref{Figure5}).
The Riemann mapping theorem tells us that there exists a biholomorphism:
\begin{align*}
  s : H^- \to \mathcal{T}
\end{align*}
and that any such map extends to a homeomorphism between the closure of these domains
and so identifies the boundaries 
$\partial H^- = \R \cup \{\infty\}$ and $\partial\mathcal{T}$.
Furthermore since $\partial\mathcal{T}$ is formed by arcs of circles (or line
segments), the extended map is regular except at the three points of $\R \cup \{\infty\}$
which correspond to the non-smooth points of $\partial\mathcal{T}$, i.e., to
the vertices $0, t, h$. The key is to understand the behaviour of $s$ at these
singular points.

First we fix the locations of the singular points. The group of holomorphic
automorphisms $\Aut(H^-)$ has a natural action on $\partial H^-$ and any element
is uniquely determined
by its images of the points $0, 1, \infty \in \partial H^-$ which can
be any three (distinct) points. Since any two biholomorphisms $H^- \to \mathcal{T}$
are related by an element of $\Aut(H^-)$, we can place the singular points of $s$
anywhere on $\partial H^-$ and once we have done this, $s$ is
uniquely specified. We place the singular points corresponding to
$h, t, 0 \in \partial \mathcal{T}$ at $0, 1, \infty$ respectively.
By uniqueness $s$ must be the inverse for $I$ restricted to $\mathcal{T}$.

We now discuss the behaviour of $s$ at the singular points. Since the corresponding
points on $\mathcal{T}$ are intersections of arcs of circles and are thus conformal
to intersections of straight lines meeting at the same angles, the behaviour
of $s$ is necessarily of the form\footnote{See \cite{MR0045823} for the details (it
is a nice application of the Schwarz reflection principle).}:
\begin{alignat}{3}\label{local_singular_behaviour}
  s &  = Z^{1/\nu_2}s_2     && \mbox{\quad near 0}\nonumber\\
  s &  = (1-Z)^{1/\nu_1}s_1 && \mbox{\quad near 1}\\
  s &  = Z^{-1/\nu_3}s_3    && \mbox{\quad near $\infty$}\nonumber
\end{alignat}
for local functions $s_i$ which are regular and non-vanishing
at the corresponding singular points.

We now consider the Schwarzian derivative of $s$. Like $s$, it is regular on $H^- \cup \partial H^-$
except possibly at the singular points $0, 1, \infty$. Calculating $\Schwarzian s$ using
the local models \eqref{local_singular_behaviour} we find that:
\begin{alignat}{2}\label{schwarzian_behaviour}
  \Schwarzian s - \frac{1 - 1/\nu_2^2}{2Z^2} - \frac{\beta_0}{Z}       & \mbox{\quad is regular at 0}\nonumber\\
  \Schwarzian s - \frac{1 - 1/\nu_1^2}{2(1-Z)^2} - \frac{\beta_1}{1-Z} & \mbox{\quad is regular at 1}\nonumber\\
  \Schwarzian s - \frac{1 - 1/\nu_3^3}{2Z^2}                & \mbox{\quad is regular at $\infty$ and}\nonumber\\
                                                 & \mbox{\quad has a zero of order 3 there}
\end{alignat}
for real constants $\beta_0, \beta_1$. In particular $\Schwarzian s$ is regular at $\infty$ and
so:
\begin{align*}
  \Schwarzian s - \frac{1 - 1/\nu_2^2}{2Z^2} - \frac{\beta_0}{Z} -
                   \frac{1 - 1/\nu_1^2}{2(1-Z)^2} - \frac{\beta_1}{1-Z}
\end{align*}
is regular at $0, 1, \infty$ and so on $H^- \cup \partial H^-$. Since it is real-valued
on $\partial H^-$, it must be constant (by Schwarz reflection).
The only way this can be compatible with the existence of triple zero noted in
\eqref{schwarzian_behaviour} is if this constant is zero and:
\begin{align*}
  \beta_0 = \beta_1 = \frac{1 - 1/\nu_1^2}{2} + \frac{1 - 1/\nu_2^2}{2} - \frac{1 - 1/\nu_3^2}{2}
\end{align*}
Using these values, we thus obtain the desired differential equation for $s$:
\begin{align}\label{schwarzian_diff_eqn}
  \Schwarzian s = \frac{1 - 1/\nu_1^2}{2(1-Z)^2} + \frac{1 - 1/\nu_2^2}{2Z^2} +
       \frac{1 - 1/\nu_1^2 - 1/\nu_2^2 + 1/\nu_3^2}{2Z(1-Z)}
\end{align}

In general, solutions to the differential equation $\Schwarzian g = h$ may be obtained as a ratio
of linearly independent solutions to associated second-order ODEs. In our
case, an elementary computation reveals that a ratio of linearly independent
solutions to the hypergeometric differential equation:
\begin{align}\label{HG_diff_eqn}
  Z(1-Z)f'' + \left(c - (a + b + 1)Z\right)f' - abf = 0
\end{align}
solves \eqref{schwarzian_diff_eqn} iff:
\begin{center}
  \begin{tabular}{ccc}
    $c - a - b = \pm 1/\nu_1$ & $1-c = \pm 1/\nu_2$ & $a - b = \pm 1/\nu_3$
  \end{tabular}
\end{center}
and furthermore all solutions may be obtained this way since $\Schwarzian g_1 =
\Schwarzian g_2$ iff $g_1, g_2$ are related by a M\"obius transformation.
For the sake of definiteness, we will take the values of $a, b, c$ given by
using the $+$ signs in the above three equations. In other words, we take:
\begin{center}
  \begin{tabular}{ccc}
    $a = \frac{1}{2}\left(1 - \frac{1}{\nu_1} - \frac{1}{\nu_2} + \frac{1}{\nu_3}\right)$&
    $b = \frac{1}{2}\left(1 - \frac{1}{\nu_1} - \frac{1}{\nu_2} - \frac{1}{\nu_3}\right)$&
    $c = 1 - \frac{1}{\nu_2}$
  \end{tabular}
\end{center}
Now \eqref{HG_diff_eqn} has regular singular points at $0, 1, \infty$ and there is
a natural basis of solutions associated to each regular singular point, obtained
by employing the method of Frobenius. We shall use the basis associated to $\infty$.
As seen by elementary computation, this basis is:
\begin{align*}
  v_1(Z) &= Z^{-a} \HGF(a, 1+a-c; 1+a-b; Z^{-1})\\
  v_2(Z) &= Z^{-b} \HGF(b, 1+b-c; 1+b-a; Z^{-1})
\end{align*}
where $\HGF$ is Gauss's hypergeometric series:
\begin{align*}
  \HGF(a, b; c; Z) &= 1 + \sum_{n\ge 1} \frac{(a)_n (b)_n}{(c)_n} \frac{Z^n}{n!}
\end{align*}
and $(q)_n = q(q+1) \cdots (q+n-1)$.

This series, with radius of convergence 1, has an analytic continuation to
the complement of any path joining two regular singular points; the standard choice,
which we follow, is to use the continuation to $\C - [1, \infty)$. In fact it is
easy to see how this works: the Frobenius method allows us to find the 
bases of solutions of \eqref{HG_diff_eqn} associated to $0, 1$ and these can be
expressed in terms of the series $\HGF$ with arguments $Z, 1-Z$ respectively.
Since the circles of convergence for the bases associated to $0, \infty$ both meet the
circle of convergence for the basis associated to $1$,
there must be a linear combination of $\HGF$ in terms of the bases elements
associated to $1$ and from there to those associated to $\infty$. The coefficients
which appear in these linear relationships are known as Kummer's connection formulae.
If, by a slight abuse of notation, we use the same symbol $\HGF$ to denote
the analytic continuation then we can present the key Kummer connection formula:
\begin{align*}
  \HGF(a, b; c; Z) = \frac{\Gamma(c)\Gamma(b-a)}{\Gamma(b)\Gamma(c-a)} (-1)^{-a} v_1(Z) +
                    \frac{\Gamma(c)\Gamma(a-b)}{\Gamma(a)\Gamma(c-b)} (-1)^{-b} v_2(Z) 
\end{align*}
For further details we recommend\footnote{We should note that although
\cite{WhittakerWatson} contains a good and thorough account, it does contain
some unfortunate sign errors.} \cite{WhittakerWatson}.

To finish, we show that the map we seek is:
\begin{align}\label{icos_inverse}
  s(Z) &= 1728^{-1/5} \frac{v_1(Z)}{v_2(Z)}
       = \frac{\HGF(\frac{11}{60}, \frac{31}{60}; \frac{6}{5}; Z^{-1})}
               {(1728Z)^{1/5}\HGF(-\frac{1}{60}, \frac{19}{60}; \frac{4}{5}; Z^{-1})}
\end{align}
(where $Z^{1/5}$ is defined using the principal branch of $\log$ on $\C - (-\infty, 0]$).
We know that the map we seek is a M\"obius transformation of $s$:
\begin{align*}
  \frac{\alpha s + \beta}{\gamma s + \delta}
\end{align*}
While $s$ is not regular at $0, 1, \infty$ its value does exist at these points and
we could use Kummer's connection formulae to verify that \eqref{icos_inverse} sends
these points to the appropriate vertices. However this is a rather involved calculation
(involving non-trivial $\Gamma$-function identities) and so we proceed differently.
First note that since $s(\infty) = 0$ we must have $\beta = 0$ and $\alpha \ne 0$.
We can thus assume $\alpha = 1$. To determine $\gamma, \delta$,
let $z_1, z_2$ be the numerator, denominator respectively in
\eqref{icos_inverse} and substitute into the identity:
\begin{align*}
  H^3(z_1, \gamma z_1 + \delta z_2) = 1728Z\cdot f^5(z_1, \gamma z_1 + \delta z_2)
\end{align*}
expanding the series $\HGF$ in \eqref{icos_inverse} to order $Z^{-1}$.
Comparison of leading terms yields $\gamma = 0, \delta^5 = 1$. Finally note that
any such value of $\delta$ will provide an inverse for $I$ since multiplication
by $e^{2\pi i/5}$ is an icosahedral rotation\footnote{For the especially
dedicated reader who desires not just an inverse but to know that \eqref{icos_inverse}
really is the inverse mapping to $\mathcal{T}$ when $\delta = 1$, the easiest
way to show this seems to be to use the reality
of the $\Gamma$-functions appearing in the Kummer connection formulae joining bases associated to
$0$ and $\infty$.}.

Using a similar expression for $\im(Z) > 0$, we could extend this function
to the open set $H^+ \cup H^- \cup (0, 1)$ so that we would have an inverse
for the restriction of $I$ to the interior of the triangle with vertices $0, \epsilon+
\epsilon^{-1}, h$.

\section{Further properties and parting words}
Our focus in these notes has been to present the icosahedral
solution of the quintic as concisely as possible, subject to the conditions
of remaining as explicit as \cite{Klein} and as self-contained as possible. As a result we have been
forced us to omit discussion of many related matters. We comment briefly on some of these
here (working over $\C$).

\subsection{Bring's curve and Kepler's great dodecahedron}\label{BringCurveSection}
  We mentioned in section \ref{Tschirnhaus_sect} that it is possible
  to reduce the general quintic to the so-called Bring-Jerrard form:
  \begin{align*}
    y^5 + y + \gamma = 0
  \end{align*}
  but that we would work with the quintic in the form \eqref{CanonicalQuintic}. We did
  this because we were following \cite{Klein}, because \eqref{CanonicalQuintic} is more general
  and because it is easy to bring out the icosahedral
  connection using the $A_5$ actions on the lines in the doubly-ruled quadric surface.
  However there is an appealing way to connect the icosahedron
  with the quintic in Bring-Jerrard form which is worth mentioning.
  The construction below is described in \cite{MR511743}.

  Firstly note that the family of quintics in Bring-Jerrard form is the smooth genus 4 curve $B$
  cut out of $\PP^4$ by the equations $\sum y_i = \sum y_i^2 = \sum y_i^3 = 0$. This is
  known as the Bring curve and has automorphism group $S_5$ corresponding to the general Galois
  group. The branched covering $B \to B/A_5 \simeq \PP^1$ allows us to define an
  invariant as before.

  Secondly, starting with an icosahedron in $\R^3$ we form Kepler's great dodecahedron $G_D$. This
  regular
  solid, which self-intersects in $\R^3$, has one face for each vertex of the icosahedron.
  It is formed by spanning the five neighbouring vertices of each vertex of the icosahedron with
  a regular pentagon and then dismissing the original icosahedron. $G_D$ thus has the same 12 vertices
  and 30 edges as the icosahedron but only 12 faces. Projection onto the common
  circumsphere $S$ yields a triple covering $G_D \to S$ with a double branching at the 12 vertices and
  after identifying $S$ with $\PP^1$ provides $G_D$ with a complex structure. Evidently $G_D$ has
  Euler characteristic $-6$ and so genus 4. In fact, as explained in \cite{MR511743}, $G_D$ is
  isomorphic to the Bring curve.

  The isomorphism $G_D \simeq B$ can be used to bring out the relationship
  between the quintic and the icosahedron.

\subsection{Modular curves and Ramanujan's continued fraction}\label{ModularSect}
  From one point of view, the exceptional geometry of the quintic is a result of the
  exceptional isomorphism:
  \begin{align*}
    A_5 \simeq PSL(2, 5)
  \end{align*}

  Corresponding to the exact sequence defining the level-5 principal congruence
  subgroup of the modular group:
  \begin{align*}
    0 \to \Gamma(5) \to PSL(2, \Z) \to PSL(2, 5) \to 0
  \end{align*}
  there is a factorization of the modular quotient:
  \begin{align*}
    j : H^* \overset{j_5} \longrightarrow X(5) \overset{\hat I}\longrightarrow X(1)
  \end{align*}
  where $H^* = H^+ \cup \Q\PP^1$ is the upper half-plane together with the $PSL(2, \Z)$ orbit
  of $\infty$ and $X(N)$ is the compactified modular curve of level $N$. The curves $X(5), X(1)$
  are rational and the map $\hat I : X(5) \to X(1)$ is a quotient by $PSL(2, 5)$ and is thus
  an icosahedral quotient. We can use this to find an inverse for the icosahedral function $I$
  in terms of Jacobi $\vartheta$-functions (provided we are willing to invert Klein's
  $j$-invariant). Indeed the map $j_5$ may be expressed as\footnote{Those comparing with
  \cite{Klein} should note that Klein's version of \eqref{inv_using_jtheta} contains some typos.}:
  \begin{align}\label{inv_using_jtheta}
    j_5(\tau) = q^{2/5}\frac{\sum_{\Z}q^{5n^2+3n}}{\sum_{\Z}q^{5n^2+n}}
              = q^{-3/5}\frac{\vartheta_{1}(\pi\tau; q^5)}{\vartheta_{1}(2\pi\tau; q^5)}
  \end{align}
  where $q=e^{\pi\tau i}$ and we are using the $\vartheta$-function notational conventions
  of \cite{WhittakerWatson}. Thus given $Z$ as in section \ref{icos_inv_sect}, if $\tau$
  satisfies $j(\tau) = 1728Z$ then $z = j_5(\tau)$ is a solution to $I(z) = Z$.

  In fact there is another expression for $j_5$, it is none other than
  Ramanujan's continued fraction:
  \begin{align*}
    j_5(\tau) = \cfrac{q^{1/5}}{1+\cfrac{q}{1+\cfrac{q^2}{1+\cfrac{q^3}{1+\cdots}}}}
  \end{align*}
  Furthermore because we know that the icosahedral vertices, edge midpoints and face centres
  in $X(5)$ lie above the points $\infty, 1, 0$ in $X(1)$, we can calculate the values
  of this continued fraction at those orbits in $H^*$ which $\frac{1}{1728}j$ maps to $\infty, 1, 0$.
  For example $j(i) = 1728$ and the corresponding edge midpoint
  equality:
  \begin{align*}
    j_5(i) = t = \sqrt{\frac{5+\sqrt{5}}{2}} - \frac{1+\sqrt{5}}{2}
  \end{align*}
  is one of the identities that famously caught Hardy's eye when Ramanujan
  first wrote to him. A beautiful account of these results
  together with a proof of \eqref{inv_using_jtheta} can be found in \cite{MR2133308}.

\subsection{Parting words}
  There is of course much more to say beyond even those remarks in sections
  \ref{BringCurveSection} and \ref{ModularSect} above.
  For example:
  \begin{itemize}
    \item There is a beautiful algorithm for solving the quintic based on iterating
          a rational function with icosahedral symmetry discovered by Doyle and
          McMullen \cite{MR1032073}.
    \item The rational parameterization of the singularity $T^2 + H^3 = 1728f^5$ we have described
          can be used to find solutions of the Diophantine equation $a^2 + b^3 + c^5 = 0$.
          See Beukers \cite{MR1487980} for details.
    \item The icosahedral solution of the quintic
          is not usually the most efficient technique for finding the roots. More practical
          formulae appear in \cite{MR1731613} for example.
  \end{itemize}

\appendix
\section{An earlier solution}\label{eisenstein_appendix}
In addition to the techniques described above, there is another approach to the
solution of the quintic discovered by Lambert\footnote{It should be pointed out that
Lambert would not have been aware that his method provided a solution of the general
quintic since the reduction the Bring-Jerrard form was not known in his time (nor was
the non-existence of a radical solution known).}
\cite{Lambert} in 1758 and again by Eisenstein \cite{Eisenstein} in 1844.

Consider the quintic in Bring-Jerrard form (up to a sign):
\begin{align}\label{BJ_form_PS}
  y^5 - y + \gamma &= 0
\end{align}
Viewing $y$ as an analytic function of $\gamma \in \C$, we claim that the branch of $y$ such
that $y(0) = 0$ has power series:
\begin{align}\label{QuinticPowerSeries}
  y(\gamma) = \sum_{k\ge 0} \binom{5k}{k}\frac{\gamma^{4k+1}}{4k+1}
\end{align}
(This can also be expressed in terms of
the generalized hypergeometric series as: $y(\gamma) =
{}_4F_3\left(\frac{4}{5},\frac{3}{5},\frac{2}{5},\frac{1}{5};
\frac{5}{4},\frac{3}{4},\frac{1}{2};5(\frac{5\gamma}{4})^4\right)\gamma$.)

This appealing result is established using analytic methods (Lagrange inversion)
in \cite{MR1057021} and \cite{MR1336074} as well as \cite{MOLagrangeInversion}.
However since this statement is really an identity of binomial coefficients it is desirable
to have a combinatorial proof for the identity \eqref{BJ_form_PS}
satisfied by the generating function \eqref{QuinticPowerSeries}.

Now the coefficients in \eqref{QuinticPowerSeries} are a special case of the sequence:
\begin{align*}
  {}_pd_k = \frac{1}{(p-1)k+1}\binom{pk}{k}
\end{align*}
which specializes to the Catalan numbers for $p=2$.
This sequence, considered long ago by Fuss \cite{Fuss}, was studied in some detail in
\cite{MR1098222}.
Just as various identities for the Catalan numbers can be established
by observing that the they count (amongst many other things) certain
lattice paths, so too can
those identities for ${}_pd_k$ which we seek for $p=5$ be established by demonstrating
that these coefficients count certain paths introduced in \cite{MR1098222}.

Although the results in \cite{MR1098222} thus provide a combinatorial proof of the generating
function identity \eqref{BJ_form_PS}, there is a more direct combinatorial  proof
presented in \cite{MR1397498} based on an observation of Raney \cite{MR0114765}. He noticed
that if $a_1, \ldots, a_m$ is any sequence of integers that sum to 1, then exactly
one of the $m$ cyclic permutations of this sequence has all of its partial sums positive.
With this in mind we consider the problem of counting sequences $a_0, \ldots, a_{kp}$
such that:
\begin{itemize}
  \item $a_0 + \cdots + a_{kp} = 1$
  \item All partial sums are positive
  \item Each $a_i$ is either 1 or $1-p$
\end{itemize}
Using Raney's observation it is clear that the number
of such sequences is ${}_pd_k$. The natural recursive structure of such sequences
provided by concatenation of $p$ such sequences, followed by a terminating value of $1-p$
then corresponds to the identity we seek. The interested reader will find details in
\cite{MR1397498}.

\bibliographystyle{plain}
\bibliography{notes}

\end{document}